\documentclass[a4paper,11pt]{article}

\usepackage[T1]{fontenc}
\usepackage[english]{babel}
\usepackage{amsmath} 
\usepackage{amssymb} 
\usepackage{amsthm} 
\usepackage{amsfonts} 
\usepackage{graphicx} 
\usepackage{latexsym} 
\usepackage{epsf} 
\usepackage{fancyhdr} 
\makeatletter

\theoremstyle{definition}    
\newtheorem{thm}{Theorem}[section]
\numberwithin{equation}{section} 
\numberwithin{figure}{section} 
\theoremstyle{plain}    
\theoremstyle{remark}
\newtheorem{rem}[thm]{\bf Remark}
\newtheorem{prop}[thm]{\bf Proposition}
\newtheorem{cor}[thm]{\bf Corollary}
\makeatother

\newcommand{\Chi}{{\mathcal X}}
\newcommand{\V}{{\mathcal V}}
\newcommand{\F}{{\mathcal F}}

\title{Normal forms of Poisson structures near a symplectic leaf}

\author{Olivier BRAHIC.\thanks{GTA, D\'epartement de Math\'ematiques, case 51,
Universit\'e Montpellier II, Place Eug\`ene Bataillon, 34095
Montpellier cedex 5, France. e-mail: brahic@math.univ-montp2.fr}}
\begin{document}
\maketitle
\begin{abstract}
In this paper, we show how one can handle the formalism developped by Yurii Vorobjev in order to give general results about the problems of linearisation and of normal form of a Poisson structure in the neighborhood of one of its symplectic leaves.
\end{abstract}

\section{Coupling tensors}
Here, we summarize the tools constructed by Vorobjev \cite{Vo} that will be useful for our purpose. 
In order to study a Poisson structure $\Pi$ on the neighborhood of one of its compact leaves $S$, we can choose a tubular 
neighborhood $p:E\longmapsto S$, seen as a vector bundle, and consider that $\Pi$ is defined on $E$. Denote $n$ the rank 
of $E$, and $2s$ the dimension of $S$, while $\Chi(S)$ and $\Omega(S)$ denote respectively the $C^\infty(S)-$modules of 
vector fields and of $1-$forms on $S$. Let $Ver:=\ker p_\star\subset TE$ denote the sub-bundle of vertical vector fields, 
and $Ver^0 \subset T^\star E$ its annihilator.

As  $\Pi$ is non-degenerate on $S$, one can always choose $E$ such that that it is {\it horizontally non-degenerate}, 
that is, verifies:
 \begin{eqnarray} \label{nondeg} \left\{ \begin{array}{ccc}
    \Pi^\sharp(Ver^0)\cap Ver    &=& \{0\} \\
    \mbox{dim }\Pi^\sharp(Ver^0) &=& \mbox{dim }S  \\
              Im\Pi^\sharp_{|S}  &=& TS,
 \end{array} \right. \end{eqnarray}
where, as usual, $\Pi^\sharp:T^\star E\longmapsto TE$ denotes the natural contraction $\alpha\longmapsto i_\alpha \Pi$.  
By putting $Hor:=\Pi^\sharp(Ver)$ one defines a subbundle of $TE$ complementary to $Ver$, or equivalently, an {\it Ehresmann connection}, 
that is $\Gamma \in \Omega(E)\bigotimes_{C^\infty(S)}  Ver$ such that $\Gamma(X_V)=X_V$ for all $X_V \in Ver$, $\Gamma$ being just the projection onto $Ver$ with kernel $Hor$. For any Ehresmann connection, given $X\in\Chi(S)$, there exists a unique section of $Hor$ that 
$p-$projects onto $X$, one calls it the {\it horizontal lift} of $X$ and denotes it $hor (X)$. The curvature of $\Gamma$ is an element of $\Omega^2(S)\otimes_{C^\infty(S)}Ver$, 
defined as $Curv_\Gamma(u,v):=[hor u,hor v]-hor[u,v]$. Also, the covariant derivative attached to $\Gamma$, 
$\partial_\Gamma: \Omega^k(S)\bigotimes_{C^\infty(S)}  C^\infty(E) \longmapsto \Omega^{k+1}(S)\bigotimes_{C^\infty(S)}  C^{\infty}(E)$ defined by the following:
$$ \begin{array}{r}
(\partial_\Gamma \F)(u_0,u_1,\dots,u_k) := \sum_{i=0}^k (-1)^i {\mathcal L}_{hor(u_i)} \F (u_0,u_1,\dots,\hat{u}_i,\dots,u_k) \\
 +\sum_{0\leq i<j\leq k} (-1)^{i+j} \F ([u_i,u_j],u_0,u_1,\dots,\hat{u}_i,\dots,\hat{u}_j,\dots,u_k).
 \end{array}$$
 
  {\it Geometric data} on $E$ are defined to be triples $(\Gamma,\V,\F)$ with $\Gamma$ an Ehresmann connection, $\V\in\Lambda^2Ver$
 a vertical $2$-vector field on $E$, and $\F\in \Omega^2(S)\bigotimes_{C^\infty(S)}  C^\infty(E)$ a {\it non-degenerate} $2$-form on $S$ 
 with values in $C^\infty(E)$.
 The first result of Vorobjev consists into describing any horizontally non-degenerate $2-$vector field on $E$ as a geometric data.
 \begin{prop} {\it
 Geometric data on $E$ are equivalent to horizontally non-degenerate $2$-vector fields}
 \end{prop}
 {\bf proof:} Given a non-degenerate $2-$vector field $\Pi$ on $E$ (not supposed to be Poisson), we define $\Gamma$ to be the 
 Ehresmann connection with kernel $Hor:=\Pi^\sharp(Ver)$, as explained above. Let 
 $\Gamma^{\Lambda2}:\Lambda^2TE\longmapsto Ver$ the external second power of $\Gamma$, define $\V\in\Lambda^2 Ver$ 
 as $\V:=\Gamma^2(\Pi)$, there only remains to define $\F$.
 
 As $\Pi$ is horizontally non-degenerate, $\Phi:=\Pi^\sharp_{|Ver^0}\longmapsto Hor$ is an isomorphism, so that one can define  $\F\in \Omega^2(S)\bigotimes_{C^\infty(S)}  C^\infty(E)$ by requiring:
 $$ \F(u,v)=-<\Phi^{-1}(hor u),hor v>, \quad \forall u,v \in \Chi(S).$$
  
 We refer the interested reader to \cite{Vo} for the construction of the non-degenerate vector field given a geometric 
 data $(\Gamma,\V,\F)$, prefering to describe the situation locally. $\Box$
 
 Let   $(x_1\dots x_{2s})$ be a coordinate chart on 
 an open neighborhood $U$ of $S$ such that $p^{-1}(U)=U\times{\mathbb R}^n$, with $(y_1\dots y_n)$ linear coordinates on 
 the fibers. We can locally express $\Pi$ as:
$$ \Pi=\sum_{i,j=1\dots2s} \Pi^X_{i,j}\partial x_i\wedge\partial x_j+\sum_{i=1\dots 2s,k=1\dots n} 2\Pi^{XY}_{i,k}\partial x_i\wedge\partial y_k+\sum_{k,l=1\dots n} \Pi_{k,l}^{Y} \partial y_k\wedge \partial y_l,$$
with $\Pi^X_{i,j}$, $\Pi_{i,k}^{XY}$, $\Pi_{k,l}^{Y}$ $\in C^\infty (p^{-1}(U))$ satisfying $\Pi^X_{i,j}=-\Pi_{j,i}^X$, and $ \Pi^Y_{k,l}=-\Pi_{l,k}^Y$.
The horizontal non-degeneracy condition means that the matrix $(\Pi_{i,j}^X)$ is invertible, let $({\Pi_{i,j}^X}^{-1})$ denote its inverse, we have:
$$\begin{array}{ccl}
 Hor&=&Span\mbox{ }\bigl\{\partial x_i+\sum_{j=1\dots 2s,k=1\dots n} {\Pi^X_{i,j}}^{-1}\Pi^{XY}_{j,k}\partial y_k, \mbox{ }i=1\dots n \bigr\}\\
    &=&Span\mbox{ }\bigl\{\partial x_i+\sum_{k=1\dots n} \beta_{i,k}\partial y_k, \mbox{ }i=1\dots n \bigr\},\\
\end{array}$$
where $\beta_{i,k}:=\sum_{j=1\dots2s} {\Pi^X_{i,j}}^{-1}\Pi_{j,k}^{XY},$
so that we get the following local formula for horizontal lifts:
$$ hor(\partial x_i)=\partial x_i+\sum_{k=1\dots n} \beta_{i,k}\partial y_k .$$
Let us denote horizontal lifts $X_i:=hor(\partial x_i)$ for $i=1\dots 2s.$
Putting  $$ {\mathcal V}_{p,q}:=\Pi_{p,q}^Y-\frac{1}{2}\sum_{l=1\dots2s} \beta_{l,p} \Pi^{XY}_{l,q}-\beta_{l,q}\Pi^{XY}_{l,p},$$
allows us to include mixed terms $\Pi^{XY}_{i,k}\partial x_i\wedge\partial y_k$ into basic ones $\Pi^X_{i,j}\partial x_i\wedge\partial x_j$:
$$\Pi=\sum_{i,j=1\dots 2s} \Pi_{i,j}^X X_i\wedge X_j + \sum_{k,l=1\dots n} {\mathcal V}_{k,l} \partial y_k \wedge \partial y_l.$$
This way, ${\mathcal V}_{p,q}$ appear to be the local coefficients for $\V$, while $\F$ expresses as:
$${\mathcal F}(\partial x_i,\partial x_j)=-\frac{1}{2}{\Pi_{i,j}^X}^{-1}.$$
 Note that the above decomposition 
 \begin{eqnarray}\Pi=\Pi_H+\V \label{decompositionhorver}\end{eqnarray} with $\Pi_H=\sum_{i,j=1\dots 2s} \Pi_{i,j}^X X_i\wedge X_j\in\Lambda^{2}Hor$ is {\it global}
  on $E$.
 
 Next proposition gives the Poisson condition in terms of geometric data. 
\begin{prop}\label{poisscond}
  {\it   Let $\Pi$ be a horizontally non-degenerate 2-vector field on $E$ with corresponding geometric data $(\Gamma,{\mathcal V},\F)$. 
     Then $\Pi$ is Poisson if and only if
$$\begin{array}{ccllc}
           [{\mathcal V},{\mathcal V}] &=&0,                                                   & &(i)\\
     \mathcal{L}_{hor(X)} {\mathcal V} &=&0 &    \forall X \in{\mathcal X}(S),\quad             &(ii)\\
           \partial_\Gamma {\mathcal F}&=&0,&                                                  &(ii)\\        
                       Curv_\Gamma(u,v)&=&{\mathcal V}^\sharp(d{\mathcal F}(u,v))&\forall u,v \in \Chi (S).\quad\quad\quad &(iv)            \\
\end{array}$$ }
\end{prop}

The proof can be carried out with a long calculus, so I refer mistrustful (and courageous) readers to my Phd. Thesis \cite{Br}.

\begin{prop}{\bf(Semi-local splitting)} \label{courburenulle}
{\it Let $\Pi$ be an horizontally non-degenerate Poisson, with corresponding geometric data 
$(\Gamma,{\mathcal V},{\mathcal F})$. In the decomposition \ref{decompositionhorver}:
$$ \Pi=\Pi_H+{\mathcal V},$$
the 2-vector field $\Pi_H$ is Poisson if and only if  $\Gamma$ has null curvature.}
\end{prop}
The proof lies in \cite{Br}, and is left to the reader.
\begin{rem}
 The problem of finding a tubular neighborhood $p:E\longmapsto S$
such that $\Pi_H$ is Poisson is the semi-local analog of Weinstein's splitting theorem (see \cite{We1},\cite{We2}), but may not hold if, for exemple, $E$ does not admit connections without curvature. However, some positive answer will be given later on. 
\end{rem} 
Before giving conditions under which we will be able to construct some homotopy between two horizontally non-degenerate Poisson structures having same vertical part, let us precise some notations. 

Given $\phi\in \Omega^1(S)\otimes_{C^{\infty}(S)}C^\infty(E)$, one defines ${\mathcal V}^\sharp(d\phi)_h \in \Omega^1(E)\otimes_{C^\infty(E)}\chi(Ver)$
 by putting
$$ {\mathcal V}^\sharp(d\phi)_h(X):={\mathcal V}_x^\sharp(d_x\phi(p_\star X))\in Ver_x\quad\forall x \in E, \quad \forall X\in T_x E.$$
It is easily verified that this definition point by point assures  $C^{\infty}(E)$-linearity.

Besides, given $\phi_1, \phi_2 \in \Omega^1(S)\otimes_{C^{\infty}(S)}C^\infty(E)$ two 1-forms on $S$ with values in $C^\infty(E)$, it is defined $\{\phi_1,\phi_2\}_{\mathcal V}$ to be the 2-form on $S$ with values in $C^\infty(E)$ (that is $\{\phi_1,\phi_2\}_{\mathcal V} \in \Omega^2(S)\otimes_{C^{\infty}(S)}C^\infty(E)$) such that for all $u_1,u_2 \in \chi(S)$:
$$ \{ \phi_1,\phi_2\}_{\mathcal V}(u_1,u_2):={\mathcal V}\bigl(d\phi_1(u_1),d\phi_2(u_2)\bigr)-{\mathcal V}\bigl(d\phi_1(u_2),d\phi_2(u_1)\bigr)  .$$

\begin{prop}\label{homotopie}
{\it Let $\Pi$ and $\Pi'$ two horizontally non-degenerate Poisson structure on some vector bundle $E$ over a compact base $S$,
with corresponding geometric datas $(\Gamma,{\mathcal V},{\mathcal F})$ et $(\Gamma',{\mathcal V}',{\mathcal F}')$.
Suppose that $\Pi$ and $\Pi'$ have same vertical part and coincide on $S$:
$$\begin{array}{rcl}
              {\mathcal V}'&=&{\mathcal V}\\  
{\mathcal F}'(u_1,u_2)_{|S}&=&{\mathcal F}(u_1,u_2)_{|S} \quad \forall u_1,u_2 \in \Chi(S).
\end{array}$$
If there exists some $1$-form with values in $C^\infty(E)$, $\phi \in \Omega^1(S) \otimes_{C^\infty(S)}C^\infty(E)$ such that:
$$\begin{array}{ccl}
        \Gamma'&=&\Gamma-{\mathcal V}^\sharp(d\phi)_h\\
  {\mathcal F}'&=&{\mathcal F}+\partial_\Gamma \phi+\frac{1}{2} \{ \phi,\phi\}_{\mathcal V},
\end{array}$$
then there exists neighborhoods ${\mathcal  U,\mathcal U'}$ of $S$ in $E$ and a diffeomorphism
$\Phi: {\mathcal U}\longmapsto {\mathcal U'}$ such that
$$\begin{array}{ccl}
 \Phi_\star \Pi&=&\Pi'\\
      \Phi_{|S}&=&Id_S.
\end{array}$$}
\end{prop}

Let' s just sketch the proof (see \cite{Vo} for the details). It consists in constructing some homotopy between $\Pi$ and $\Pi'$, that is just the flow of a time-dependant vector field $X_t$ linked to $\phi$ via the non-degeneracy condition for $\F$. The conditions on
$\phi$ are exactly the ones for $X_t$ to assure that it will tract $\Pi$ to $\Pi'$.

\section{Linearisation near a symplectic leaf}
The problem of linearising a Poisson structure that vanishes at one point was studied in 
\cite{Co}, \cite{Du1}, \cite{Du2}, \cite{Du3}, \cite{DAffn}... In this section, we want to 
see how these results extend to a full neighborhood of a fixed symplectic leaf, in the most general situation possible.

\begin{prop}{\bf(Linearisation of the vertical part)}\label{linearisationverticale}
{\it Let $\Pi$ an horizontally non-degenerate Poisson structure on some vector bundle $p:E\longmapsto S$
over a compact base $S$, with corresponding geometric data $(\Gamma,{\mathcal V},{\mathcal F})$. Suppose that, fore some $x\in S$, the germ
of ${\mathcal V}$  at $x$ is linearisable (as a Poisson structure on $E_x$), then the germ of $\Pi$
along $S$ is equivalent to the one of a Poisson structure with associated geometric data 
$(\Gamma',{\mathcal V}^{(1)},{\mathcal F}')$, for some $\Gamma'$ and $\F'$. Here ${\mathcal V}^{(1)}$ denotes the linear part of $\V$.}
\end{prop}
{\bf proof:} The vertical part is locally trivial (as $\V$ is invariant by a connection according to condition $(ii)$) so the proposition is trivial
if $S$ is some ball. For more general $S$, the problem is to paste local linearising charts. The argument is the following: 
instead of linearising directly $\V$, one constructs a vector field $Z\in\Chi(E)$ such that:
\begin{eqnarray}
 [Z,{\mathcal V}]=-{\mathcal V} \label{eqhom}\\
          Z^{(1)}=L, \label{lioucond}
\end{eqnarray}
where $Z^{(1)}$ denotes the part of order one in the fibers of $Z$ , and $L=\sum_{i=1\dots n} y_i\partial y_i$ the so called 
{\it Liouville vector field} (note that the form $\sum_{i=1\dots n} y_i\partial y_i$ does not depend on the chosen trivialisation of 
the vector bundle $E$ as it is equivalent to requiring that the flow $\phi_t$ of $Z^{(1)}$ is multiplication by $\exp(t)$). The vector field
$Z$ turns out to be {\it globally} linearisable, so that, up to a diffeomorphism defined on a some neighborhood of $S$, we have
 $$[L,\V]=-\V.$$
 It is then a basic fact that this last condition forces $\V$ itself to be fiberwise linear.
 
 There remains to construct $Z$: let $(U_i)_{i\in I}$ some finite cover of $S$ with $U_i$  open balls of $S$, 
that admit trivialisations  $\psi_i:p^{-1}(U_i)\longmapsto {\mathcal U}_i\times {\mathbb R}^n
$
of the vector bundle $E$, and $(\rho_i)_{i\in I}$ some subordinate partition of unity. The local triviality of $\V$, 
both with the linearisability condition enables to construct vector fields $Z_i$ over $p^{-1}(U_i)$ such that
$$\begin{array}{rcl}
 [Z_i,{\mathcal V}]&=&-{\mathcal V} \\
          Z_i^{(1)}&=&L.
\end{array}$$
Let $Z:=\sum_{i\in I} \rho_i Z_i$, then \ref{lioucond} is easily verified, and \ref{eqhom} holds because $\rho_i$ are basic and
$\V$ vertical. $\Box$
\begin{prop}{\bf(Changing the connection)}\label{changementconnectionh1}
{\it Let $\Pi$ an horizontally non-degenerate Poisson structure on some vector bundle $p:E\longmapsto S$
over a compact base $S$, with corresponding geometric data $(\Gamma,{\mathcal V},{\mathcal F})$. Suppose that
the first Poisson cohomology space $H^1({\mathcal V}_x)$ of ${\mathcal V}$ at some $x\in S$
is trivial, and let $\Gamma'$ some connection leaving ${\mathcal V}$ invariant.

Then the germ of $\Pi$ along $S$ is equivalent to the one of a Poisson structure with associated geometric data 
 $(\Gamma',{\mathcal V},{\mathcal F}')$, for some ${\mathcal F}'$.}
\end{prop}
{\bf proof:} According to \ref{homotopie}, we only have to construct some $\phi\in \Omega^1(S)\otimes_{C^\infty(S)} C^\infty(E)$ such that:
     $$ \Gamma=\Gamma'+{{\mathcal V}}^\sharp(d\phi)^h.$$
Keeping the same notations as before, over $p^{-1}(U_j)$ $\Gamma$ and $\Gamma'$ take the form 
$$\begin{array}{ccc} 
\Gamma(\partial x_i)&=&-\sum_{k=1}^{n} \beta_{i,k} \partial y_k\\
\Gamma'(\partial x_i)&=&-\sum_{k=1}^{n} \beta'_{i,k} \partial y_k.    
\end{array}$$
As $\Gamma$ and $\Gamma'$ are supposed to leave $\V$ invariant, we have 
$$\bigl[\sum_{k=1}^{n} (\beta_{i,k}-\beta'_{i,k})\partial y_k,{\mathcal V} \bigr]=0, \quad \forall i=1\dots 2s.$$
which means that $\sum_{k=1}^{n} (\beta_{i,k}-\beta'_{i,k})\partial y_k$ is a  $1$-cocycle for the Poisson cohomology
of ${\mathcal V}$ (with parameters in $U_i$). By hypothesis, there exist $\phi_i^j\in C^\infty(p^{-1}(U_j))$ such that for all $i=1\dots 2s$,  
    $$[\mathcal{V},\phi_i^j]=\sum_{k=1}^{n} (\beta_{i,k}-\beta'_{i,k})\partial y_k.$$
Let $\phi^j:= \sum_{i=1\dots 2s} \phi_i\otimes dx_i$,  one gets on $p^{-1}(U_j))$ some $1-$form such that
    $$ \Gamma=\Gamma^{(1)}+{{\mathcal V}}^\sharp(d\phi^j)^h.$$
It is easily verified that, as $\rho_i$ are basic, $\phi:=\sum_{j\in J} \rho_j \phi^j$ gives the desired 1-form.$\Box$

We are now handing enough tools to state the following theorem.

\begin{thm}{\bf(Semi-local linearisation of Poisson structures)}
{\it Let $\Pi$ an horizontally non-degenerate Poisson structure on some vector bundle $p:E\longmapsto S$
over a compact base $S$, with corresponding geometric data $(\Gamma,{\mathcal V},{\mathcal F})$. Suppose 
that the germ of ${\mathcal V}$  at $x$ is linearisable (as a Poisson structure on $E_x$) for some $x\in S$
and that the first Poisson cohomology space $H^1({\mathcal V_x^{(1)}})$ of ${\mathcal V^{(1)}}$ is trivial.

 Then the germ of $\Pi$ along $S$ is equivalent to the one of a Poisson structure with associated geometric data 
$(\Gamma^{(1)},{\mathcal V}^{(1)},{\mathcal F})$. Here ${\mathcal V}^{(1)}$ and $\Gamma^{(1)}$ respectively denote the linear parts of $\V$ and $\Gamma$.}
\end{thm}

{\bf proof:} First apply proposition \ref{linearisationverticale}: up to a diffeomorphism,  $\Pi$ has associated geometric data
$(\Gamma,{\mathcal V}^{(1)},{\mathcal F})$, then  $\Gamma$ locally takes the form 
$$ hor(\partial x_i)=\partial x_i+\sum_{j=1\dots n}  \bigl(\sum_{k=1\dots n}\beta_{i,j}^k(x)y_k+\hat{\beta}_{i,j}(x,y)\bigr) \partial y_j,$$
with $\hat{\beta}_{i,j}$ of order higher than one in the $y_i-$variables. Just isolating terms of order one in the equation
$${\mathcal L}_{hor(\partial x_i)}\V^{(1)}=0,$$ one sees that $\Gamma^{(1)}$ leaves $\V$ invariant, so we only have to apply
 proposition \ref{changementconnectionh1}.$\Box$

 We can immediately derive the following corollaries.
   
\begin{cor}{\bf(Semi-local non-degeneracy of Poisson structures)}
{\it Let $\Pi$ an horizontally non-degenerate Poisson structure on some vector bundle $p:E\longmapsto S$
over a compact base $S$, with corresponding geometric data $(\Gamma,{\mathcal V},{\mathcal F})$. Suppose 
that the vertical linear part ${\mathcal V}_x^{(1)}$  at some $x\in S$ is non-degenerate
and that the first Poisson cohomology space $H^1({\mathcal V_x^{(1)}})$ of ${\mathcal V^{(1)}}$ is trivial. 

Then the germ of $\Pi$ along $S$ is equivalent to the one of a Poisson structure with associated geometric data 
$(\Gamma^{(1)},{\mathcal V}^{(1)},{\mathcal F}')$ for some $\F'$.}
\end{cor}

\begin{cor}{\bf(Semi-local non-degeneracy of transversally compact semi-simple Poisson structures)}
{\it  Let $\Pi$ an horizontally non-degenerate Poisson structure on some vector bundle $p:E\longmapsto S$
over a compact base $S$, with corresponding geometric data $(\Gamma,{\mathcal V},{\mathcal F})$. Suppose 
that the vertical linear part ${\mathcal V}_x^{(1)}$  at some $x\in S$ is associated to a compact semi-simple
Lie algebra.
  
Then the germ of $\Pi$ along $S$ is equivalent to the one of a Poisson structure with associated geometric data 
$(\Gamma^{(1)},{\mathcal V}^{(1)},{\mathcal F}')$ for some $\F'$.}
\end{cor}

\section{Semi-local splitting}
Last theorem enables us to write certain Poisson structures $\Pi$ as 
$$ \Pi=\Pi_H+{\mathcal V}^{(1)}$$
in a neighborhood of some leaf $S$, with $\Pi_H$ inducing some linear conection on $E$. That' s what we mean by {\it global linearisation}. But up to now, we didn' t 
care about the 2-form $\F$, the first question being: what can we expect of it ? Some motivation is given by proposition
\ref{courburenulle}, as we can see in proposition \ref{poisscond} that the curvature of $\Gamma$ is intimely linked 
to $\F$. We already pointed out the necessary condition that $E$ shall admit some connection without curvature, so
this condition is not restrictive in the following proposition.
\begin{thm}{\it
Let $\Pi$ an horizontally non-degenerate Poisson structure on some vector bundle $p:E\longmapsto S$
over a compact base $S$, with corresponding geometric data $(\Gamma,{\mathcal V},{\mathcal F})$. Suppose 
that the vertical  part ${\mathcal V}_x$  at some $x\in S$ has trivial first Poisson cohomology space 
$H^1({\mathcal V}_x)$, and that there exists some Ehresmann connection $\Gamma'$ without curvature leaving $\V$ 
invariant.

Then, up to a diffeomorphism on a neighborhood of $S$, $\Pi$ decomposes into 
           $$\Pi=\Pi_H+{\mathcal V},$$ 
with $\Pi_H$ and ${\mathcal V}$ Poisson structures such that  $\Pi_H$ is regular and ${\mathcal V}$ vertical, vanishing 
on $S$.} 
\end{thm}
{\bf proof:} Just apply propositions \ref{changementconnectionh1} and \ref{courburenulle}.

\end{document}